\documentclass[reqno, 10pt]{amsart}  
\usepackage{amsmath}
\usepackage{amssymb}
\usepackage{hyperref}
\hypersetup{colorlinks=true}
\usepackage{graphicx}
\graphicspath{~/Desktop/first-order hyperbolic/}
\numberwithin{equation}{section}

\newtheorem{prop}{Proposition}[section]
\newtheorem{lemma}[prop]{Lemma}
\newtheorem{thm}[prop]{Theorem}
\newtheorem{cor}[prop]{Corollary}

\newcommand{\shrink}[1]{ {\scriptstyle {\textstyle {#1} } } }
\newcommand{\smfrac}[2]{ \shrink{ \frac{#1}{#2} } }

\newcommand{\lin}{\langle}
\newcommand{\rin}{\rangle}
\newcommand{\Sym}{\mathbb{S}^{n}}
\newcommand{\Symp}{\Sym_{\plus}} 
 
\newcommand{\realsnp}{\mathbb{R}^n_{\plus}}

\newcommand{\plus}{{\scriptscriptstyle +}}

\renewcommand{\int}{\mathrm{int}}

\newcommand{\cp}{\mathrm{CP}} 

\newcommand{\lambdamin}{\lambda_{\mathrm{min}}} 
\newcommand{\grad}{\nabla} 

\newcommand{\diamval}{\mathrm{diam}_{z}  }

\usepackage{lipsum}

\newcommand\blfootnote[1]{%
  \begingroup
  \renewcommand\thefootnote{}\footnote{#1}%
  \addtocounter{footnote}{-1}%
  \endgroup
}

\begin{document}

\newpage 
$ \textrm{~} $ \quad \vspace{-3mm}

\title[A Framework for Applying Subgradient Methods]{A Framework \\ for Applying Subgradient Methods \\ to Conic Optimization Problems \\ (\MakeLowercase{version 2})$^*$}

\begin{abstract}
A framework is presented whereby a general convex conic optimization problem is transformed into an equivalent convex optimization problem whose only constraints are linear equations and whose objective function is Lipschitz continuous.  Virtually any subgradient method can be applied to solve the equivalent problem.  Two methods are analyzed. \vspace{1mm} 
\end{abstract} \vspace{-3mm}

\author[J. Renegar]{James Renegar  }
\address{School of Operations Research and Information Engineering,
 Cornell University, Ithaca, NY, U.S.}
 
\thanks{Special thanks to Yurii Nesterov for encouragement at an early, critical stage of the research. And thanks to Rob Freund, whose correspondence sparked the realization of how to best partition the results into two papers.}

\maketitle

\vspace{-7mm}

\section{{\bf  Introduction}}  \label{sect.a}
Given a conic optimization problem for which a strictly feasible point is known, we provide a  transformation to an equivalent convex optimization problem which is of the same dimension, has only linear equations as constraints (one more linear equation than the original problem), and has Lipschitz-continuous objective function defined on the whole space. Virtually any subgradient method can be applied to solve the equivalent problem, the cost per iteration dominated by computation of a subgradient and its orthogonal projection onto a subspace (the same subspace at every iteration, a situation for which preprocessing is effective). \vspace{1mm}

We develop representative complexity  results for two methods, one of which is executable under an ideal circumstance (knowing the optimal value), but the other of which is general (requiring only that a strictly feasible point be known).  \blfootnote{$^\mathbf{*}$  The development of algorithms has been streamlined and considerably strengthened.} \vspace{1mm}

Perhaps most surprising is that the transformation to an equivalent problem is simple and so is the basic theory, and yet the approach has been overlooked until now, a blind spot.  \vspace{1mm}

The following section presents the transformation and basic theory. Representative algorithmic implications  are developed in Section \ref{sect.c}.  A general example is presented in Section \ref{sect.d}, highlighting key differences with traditional literature on subgradient methods.  \vspace{1mm}

This paper significantly extends subgradient-method results first reported in \cite{renegar2014efficient} (as well as results in the previously posted version of the present paper).  A companion paper will encompass the accelerated gradient-method results reported in \cite{renegar2014efficient}.  The general theory presented in the following section is the foundation for each paper. \vspace{1mm}

\section{{\bf  Basic Theory}}  \label{sect.b}

The theory given here is elementary and yet has been overlooked, the ``blind spot'' referred to above.  \vspace{1mm}

Let $ {\mathcal E} $ be a finite-dimensional real Euclidean space.  Let $ {\mathcal K}  \subset {\mathcal E} $ be a proper, closed, convex cone with non-empty interior. \vspace{1mm}

Fix a vector $ e \in \int(  {\mathcal K} ) $ (interior).  We refer to $ e $ as the ``distinguished direction.''   For each $ x \in {\mathcal E} $, let 
\[   \lambda_{\min}(x) :=  \inf \{ \lambda: x - \lambda \, e \notin {\mathcal K}  \} \; , \]
that is, the scalar $ \lambda $ for which $ x - \lambda e $ lies in the boundary of $ {\mathcal K}  $.  (Existence and uniqueness of $ \lambda_{\min}(x) $ follows from $ e \in \int({\mathcal K} ) \neq {\mathcal E} $ and convexity of $ {\mathcal K}  $.)  \vspace{1mm}

If, for example, $ {\mathcal E} = \Sym $ ($ n \times n $ symmetric matrices), $ {\mathcal K}  = \Symp $ (cone of positive semidefinite matrices), and $ e = I $ (the identity), then $ \lambda_{\min}(X) $ is the minimum eigenvalue of $ X $. \vspace{1mm}

On the other hand, if $ {\mathcal K}  = \realsnp $ (non-negative orthant) and  $ e $ is a vector with all positive coordinates, then $ \lambda_{\min}(x) = \min_j x_j/e_j $ for $ x \in \mathbb{R}^n $. Clearly, the value of $ \lambda_{\min}(x) $ depends on the distinguished direction $ e $ (a fact the reader should keep in mind since the notation does not reflect the dependence).  \vspace{1mm}

Obviously, $ {\mathcal K}  = \{ x: \lambda_{\min}(x) \geq 0 \} $ and $ \int({\mathcal K} ) = \{ x: \lambda_{\min}(x) > 0 \} $. Also,
\begin{equation}  \label{eqn.ba}
  \lambda_{\min}(sx + te) = s \, \lambda_{\min}(x) + t \quad \textrm{for all $ x \in {\mathcal E}$  and scalars $ s \geq 0 $,  $ t $} \; .  
  \end{equation}

Let 
\[ \bar{{\mathcal B}} := \{ v \in {\mathcal E}: e + v, e - v \in {\mathcal K}  \} \; , \]
a closed, centrally-symmetric, convex set with nonempty interior. 
Define a seminorm\footnote{Recall that a seminorm $ \| \, \, \| $ satisfies $ \| tv \| = |t| \, \| v \| $ and $ \| u + v \| \leq \| u \| + \| v \| $, but unlike a norm, is allowed to satisfy $ \| v \| = 0 $ for $ v \neq 0 $.}  on $ {\mathcal E} $ according to 
\[ \| u \|_{\infty} := \min \{ t: u = tv \textrm{ for some $ v \in \bar{{\mathcal B}} $} \} \; . \]
Let $ \bar{B}_{\infty}(x,r) $ denote the closed ball centered at $ x $ and of radius $ r $. Clearly, $ \bar{B}_{\infty}(0,1) = \bar{{\mathcal B}} $, and $ \bar{B}_{\infty}(e,1) $ is the largest subset of $ {\mathcal K}  $ that has symmetry point $ e $, i.e., for each $ v $, either both points $ e + v $ and $ e - v$ are in the set, or neither point is in the set.      
\vspace{1mm}

It is straightforward to show $ \| \, \, \|_{\infty} $ is a norm if and only if $ {\mathcal K}  $ is pointed (i.e., contains no subspace other than $ \{ \vec{0}\} $). \vspace{1mm}

\begin{prop}  \label{prop.ba}
  The function $ x \mapsto \lambda_{\min}(x) $ is concave and Lipschitz continuous: 
\[    | \lambda_{\min}(x) - \lambda_{\min}(y) | \leq \| x - y \|_{\infty} \quad \textrm{for all $ x,y \in {\mathcal E} $} \; .   \]  
\end{prop} 
\noindent {\bf Proof:} Concavity follows easily from the convexity of $ {\mathcal K}  $, so we focus on establishing Lipschitz continuity.  \vspace{1mm}

Let $ x, y \in {\mathcal E} $.  According to (\ref{eqn.ba}), the difference $ \lambda_{\min}(x + te) - \lambda_{\min}(y+te) $ is independent of $ t $, and of course so is the quantity $ \| (x + te) - (y + te) \|_{\infty} \; . $ Consequently, in proving the Lipschitz continuity, we may assume $ x $ lies in the boundary of $ {\mathcal K}  $, that is, we may assume $ \lambda_{\min}(x) = 0 $. The goal, then, is to prove
\begin{equation}  \label{eqn.bb}
    | \lambda_{\min}(x + v) | \leq \| v \|_{ \infty} \quad \textrm{for all $ v \in {\mathcal E} $} \; . \end{equation} 
  
We consider two cases.  First assume $ x + v $ does not lie in the interior of $ {\mathcal K}  $, that is, assume $ \lambda_{\min}(x+v) \leq 0 $. Then, to establish (\ref{eqn.bb}), it suffices to show $ \lambda_{\min}(x + v) \geq - \| v \|_{\infty} \; , $ that is, to show
\begin{equation}  \label{eqn.bc}
    x + v + \| v \|_{\infty} \, e \in {\mathcal K}  \; . 
    \end{equation} 
However, 
\begin{equation}  \label{eqn.bd}
   v + \| v \|_{\infty} \, e \in \bar{B}_{\infty}(\| v \|_{\infty} \, e, \| v \|_{\infty}) \subseteq {\mathcal K}  \; , 
   \end{equation} 
the set containment due to $ {\mathcal K}  $ being a cone and, by construction, $ \bar{B}_{\infty}(e,1) \subseteq {\mathcal K}   $. Since $ x \in {\mathcal K}  $ (indeed, $ x $ is in the boundary of $ {\mathcal K}  $), (\ref{eqn.bc})  follows. \vspace{1mm}

Now consider the case $ x + v \in {\mathcal K}  $, i.e., $ \lambda_{\min}(x+v) \geq 0 $. To establish (\ref{eqn.bb}), it suffices to show $ \lambda_{\min}(x + v) \leq \|v\|_{\infty} \; ,  $ that is, to show
\[                x + v - \|v\|_{\infty} \, e \notin \int({\mathcal K} ) \; . \]
Assume otherwise, that is, assume 
\[    x = w + \| v \|_{\infty} \, e - v \quad \textrm{for some $ w \in \int({\mathcal K} ) $} \; . \]
Since $ \| v \|_{\infty} \, e - v \in {\mathcal K}  $ (by the set containment on the right of (\ref{eqn.bd})), it then follows that $ x \in \int({\mathcal K} ) $, a contradiction to $ x $ lying in the boundary of $ {\mathcal K}  $.    \hfill $ \Box $
 \vspace{3mm}

Assume the Euclidean space $ {\mathcal E} $ is endowed with inner product written $ u \cdot v $.  Let  $ \mathrm{Affine} \subseteq {\mathcal E} $ be an affine space, i.e., the translate of a subspace. For fixed $ c \in {\mathcal E} $, consider the conic program
\[  
  \left. \begin{array}{rl}
\inf & c \cdot x  \\
\textrm{s.t.} & x \in \mathrm{Affine}  \\
  & x \in {\mathcal K}   \end{array} \right\} \cp  \]  
Let $ z^* $ denote the optimal value. \vspace{1mm}

Assume $ c $ is not orthogonal to the subspace of which $ \mathrm{Affine} $ is a translate, since otherwise all feasible points are optimal.  This assumption implies that all optimal solutions for CP lie in the boundary of $ {\mathcal K}  $. \vspace{1mm}

Assume $ \mathrm{Affine} \cap \int({\mathcal K} ) $ -- the set of strictly feasible points -- is nonempty.  Fix a strictly feasible point, $ e $.  The point $ e $ serves as the distinguished direction. \vspace{1mm}

For scalars $ z \in \mathbb{R}  $, we introduce the affine space
\[ \mathrm{Affine}_{z} := \{ x \in \mathrm{Affine}:   c \cdot x = z \} \; . \]
Presently we show that for any choice of $ z $   satisfying $ z < c \cdot e \; ,  $ CP can be easily transformed into an equivalent optimization problem in which the only constraint is $ x \in \mathrm{Affine}_{z} \; . $ We need a simple observation.   \vspace{1mm}

\begin{lemma}  \label{lem.bb} 
Assume $ \mathrm{CP} $ has bounded optimal value.  

$ \textrm{~} $ \qquad \qquad  \qquad  If $ x \in \mathrm{Affine} $ satisfies $ c \cdot x < c \cdot e $, then $ \lambdamin(x) < 1 \; . $ 
\end{lemma}
\noindent {\bf Proof:}  If $ \lambda_{\min}(x) \geq 1 $, then $ e + t( x - e) $ is feasible for all $ t \geq 0 $ (using (\ref{eqn.ba})).  As the function $ t \mapsto c \cdot  \big( e + t(x- e) \big) $ is strictly decreasing (because $ c \cdot x < c \cdot e $), this implies CP has unbounded optimal value, contrary to assumption. \hfill $ \Box $
 \vspace{3mm}
 
For  $ x \in {\mathcal E}  $ satisfying $ \lambda_{\min}(x) < 1 $, let $ \pi(x) $ denote the point where the half-line beginning at $ e $ in direction $ x - e $ intersects the boundary of $ {\mathcal K}  $:
\[  \pi(x) :=  e + \smfrac{1}{1 - \lambda_{\min}(x)} (x - e)    \]
(to verify correctness of the expression, observe (\ref{eqn.ba}) implies $ \lambda_{\min}(\pi(x)) = 0 $). 
We refer to $ \pi(x) $ as ``the projection (from $ e $) of $ x $ to the boundary of the \vspace{1mm}  feasible region.'' 

The centrality of the following result to the development makes the result be a theorem even if the proof is straightforward. \vspace{1mm}

\begin{thm} \label{thm.bc}
Let $ z $ be any value satisfying \, $   z < c \cdot e \; . $  If $ x^* $ solves
\begin{equation}  \label{eqn.be}
  \begin{array}{rl}
   \sup & \lambdamin(x) \\
 \mathrm{s.t.} &  x \in \mathrm{Affine}_{z} \; ,  \end{array} 
    \end{equation} 
then $ \pi( x^* ) $ is optimal for $ \mathrm{CP} $. Conversely, if $ \pi^*  $ is optimal for $ \mathrm{CP} $, then $ x^* :=   e + \frac{c \cdot e - z }{c \cdot e - z^* } ( \pi^* - e) $ is optimal for (\ref{eqn.be}), and $ \pi^*  = \pi( x^* ) $.  
\end{thm}
\noindent {\bf Proof:} Fix a value satisfying $ z < c \cdot e  $. It is easily proven from the convexity of $ {\mathcal K}  $ that $ x \mapsto \pi(x) $ gives a one-to-one map from $ \mathrm{Affine}_{z} $      onto 
\begin{equation}  \label{eqn.bf}
  \{ \pi \in \mathrm{Affine} \cap  \mathrm{bdy}({\mathcal K} ) :  c \cdot \pi < c \cdot e \} \; , 
  \end{equation}
where $ \mathrm{bdy}({\mathcal K} )  $ denotes the boundary of $ {\mathcal K}  $. \vspace{1mm}

For $ x \in \mathrm{Affine}_{z} \; , $ the CP objective value of $ \pi(x) $ is 
\begin{align}
  c \cdot \pi(x) & = c \cdot \big( e+ \smfrac{1}{1 - \lambda_{ \min}(x)} (x - e) \big) \nonumber \\
                    & = c \cdot e + \smfrac{1}{1 - \lambda_{ \min }(x)} ( z - c \cdot e ) \; , \label{eqn.bg}
                    \end{align}
 a strictly-decreasing function of $ \lambda_{\min}(x) $.                    
 Since the map $ x \mapsto \pi(x) $ is a bijection between $ \mathrm{Affine}_{z} $ and the  set (\ref{eqn.bf}), the theorem readily follows. \hfill $ \Box $
 \vspace{3mm}
 
CP has been transformed into an equivalent linearly-constrained maximization problem with concave, Lipschitz-continuous objective function. Virtually any subgradient method --  rather, {\em  supgradient} method -- can be applied to this problem, the main cost per iteration being in computing a supgradient and projecting it onto the subspace $ {\mathcal L}  $ of which the affine space $ \mathrm{Affine}_{z} $ is a translate.   \vspace{1mm}

For illustration, we digress to interpret the implications of the development thus far for the linear program
\[ 
  \left. \begin{array}{rl}
 \min_{x \in \mathbb{R}^n} & c^T x \\
\textrm{s.t.} & Ax = b \\
 & x \geq 0  \end{array}   \right\} \, \mathrm{LP}  
\] 
assuming $ e = {\bf 1} $ (the vector of all ones), in which case $ \lambda_{\min}(x) = \min_j x_j \; , $ and $ \| \, \, \|_{\infty} $ is the $ \ell_{\infty} $ norm, i.e., $ \| v \|_{\infty} = \max_j |v_j| $.  Let the number of rows of $ A $ be $ m \geq 1 $. \vspace{1mm}

For any scalar $ z < c^T {\bf 1} $, Theorem~\ref{thm.bc} asserts that LP is equivalent to   
\begin{equation} \label{eqn.bh}
  \begin{array}{rl}
    \max_x & \min_j x_j \\
     \textrm{s.t.} & Ax = b \\
           & c^T x = z \; , \end{array} \end{equation} 
in that when $ x $ is feasible for (\ref{eqn.bh}), $ x $ is optimal if and only if the projection \\ $ \pi(x) = \mathbf{1}  + \smfrac{1}{1 - \min_j x_j} (x - \mathbf{1} ) $ is optimal for LP. 
  The setup is shown schematically in the following figure:

$ \textrm{~} $ \quad \qquad \quad \qquad  \qquad   \includegraphics[scale=.26]{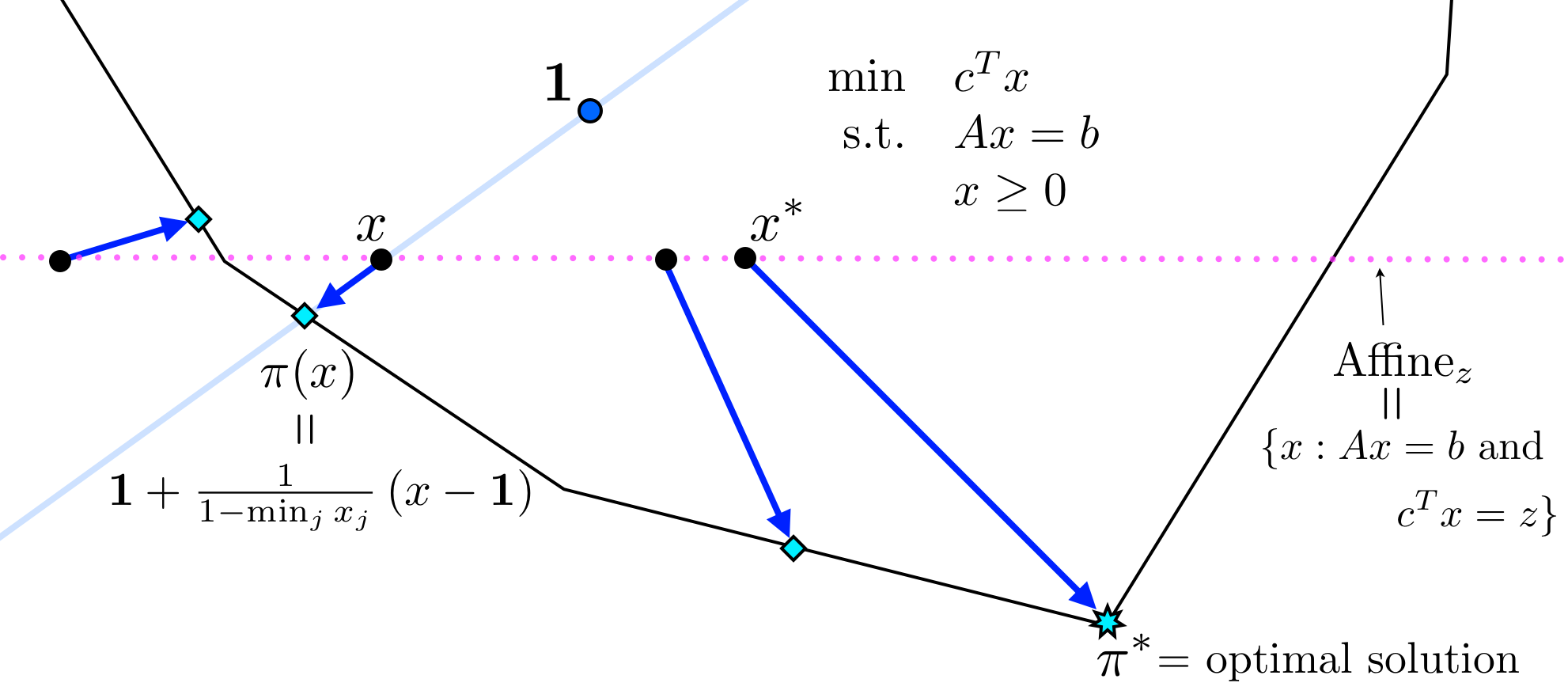}
\vspace{4mm}

Proposition~\ref{prop.ba}  asserts that, as is obviously true, $ x \mapsto \min_j x_j $ is $ \ell_{\infty} $-Lipschitz continuous with constant 1. Consequently, the function also is $ \| \, \, \|_2 $-Lipschitz continuous with constant 1, as is relevant if supgradient methods rely on the standard inner product in computing supgradients and their orthogonal projections onto the subspace $ {\mathcal L}  $    of which $ \mathrm{Affine}_{z} $ is a translate, i.e.,   $ {\mathcal L}  = \{ v: Av = 0 \textrm{ and } c^T v = 0 \} \; . $    \vspace{1mm}

With respect to the standard inner product, the supgradients of $ x \mapsto \min_j x_j $ at $ x $ are  the convex combinations of the standard basis vectors $ e(k) $ for which $ x_k = \min_j x_j $. Consequently, the projected supgradients at $ x $ are the convex combinations of the vectors $ \bar{P}_k $ for which $ x_k = \min_j x_j $, where $ \bar{P}_k $ is the $ k^{th} $ column of the matrix projecting $ \mathbb{R}^n $ onto  the nullspace of $ \bar{A} = \left[ \begin{smallmatrix}  A \\ c^T \end{smallmatrix} \right] $, that is  
\[ 
   \bar{P}  := I- \bar{A}^T (\bar{A} \,   \bar{A}^T)^{-1} \bar{A}  \; . \] 

If $ m \ll n $, then $ \bar{P} $ is not computed in its entirety, but instead the matrix $ \bar{M}  = ( \bar{A} \bar{A}^T)^{-1} $ if formed as a preprocessing step, at cost $ O(m^2 n) $. Then, for any iterate $ x $ and an index $ k $ satisfying $ x_k = \min_j x_j $, the projected  supgradient $ \bar{P}_k $ is computed according to
 \[  u = \bar{M} \, \bar{A}_k \quad \rightarrow \quad  v = \bar{A}^T u \quad \rightarrow \quad \bar{P}_k = e(k) - v \; , \]
for a cost of $ \, O(m^2 \, + \, \# \mathrm{non\_zero\_entries\_in\_} A   )  $ per iteration.  \vspace{1mm}

Before returning to the general theory, we note that if the choices are $ {\mathcal E} = \Sym $, $ {\mathcal K}  = \Symp $ and $ e = I $ (and thus $ \lambda_{\min}(X) $ is the minimum eigenvalue of $ X $), then with respect to the trace inner product, the supgradients at $ X $ for the function $ X \mapsto \lambda_{\min}(X) $ are the convex combinations of the matrices $ v v^T $, where $ Xv = \lambda_{\min}(X) v $ and $ \| v \|_2 = 1 $.   \vspace{1mm}

{\em  Assume, henceforth, that CP has at least one optimal solution, and that $ z $ is a fixed scalar satisfying \, $ z < c \cdot e $.}   Then the equivalent problem (\ref{eqn.be}) has at least one optimal solution.     Let $ x^*_{z} $ denote any of the optimal solutions for the equivalent problem, and recall
$ z^* $ denotes the optimal value of CP. A useful characterization of the optimal value for the equivalent problem is easily provided. \vspace{1mm}

\begin{lemma}  \label{lem.bd}
 \[   \lambda_{\min}(x^*_{z }) = \frac{z - z^* }{c \cdot e - 
 z^* } \] 
  \end{lemma}
\noindent {\bf Proof:}   
By Theorem~\ref{thm.bc}, $ \pi(x^*_{z}) $ is optimal for CP -- in particular, $ c \cdot  \pi(x^*_{z})  = z^* $.   Thus, according to (\ref{eqn.bg}),
\[  z^* = c \cdot e + \smfrac{1}{1 - \lambda_{\min}(x^*_{z })} \, ( z - c \cdot e) \; . \]
Rearrangement completes the proof.      \hfill $ \Box $
 \vspace{3mm}
 
We focus on the goal of computing a point $ \pi $ which is feasible for CP and has better objective value than  $ e$ in that 
\begin{equation}  \label{eqn.bi}
     \frac{c \cdot \pi - z^* }{c \cdot e - z^* } \leq \epsilon \; 
     \end{equation} 
where $ 0 < \epsilon < 1 $.  Thus, for the problem of primary interest, CP, the focus is on relative improvement in the objective value. 
\vspace{1mm}

 The following proposition provides a useful characterization of the accuracy needed in approximately solving the CP-equivalent problem (\ref{eqn.be}) so as to ensure that for the computed point $ x $, the projection $ \pi = \pi(x) $ satisfies (\ref{eqn.bi}). \vspace{1mm}

\begin{prop}  \label{prop.be}
If $ x \in \mathrm{Affine}_z $ and $ 0 < \epsilon < 1 $,  then 
\begin{align*} 
   & \frac{ c \cdot  \pi(x)  - z^* }{c \cdot e - z^* } \, \leq \,  \epsilon   
    \\ & \qquad   \qquad  \qquad   \textrm{if and only if}    \\ & \qquad  \qquad  \qquad  \qquad   
    \lambdamin( x^*_{ z}  ) - \lambdamin( x) \, \leq \, \frac{\epsilon }{1 - \epsilon } \, \, \frac{c \cdot e - z }{\, \, \, c \cdot e - z^* } \; .  
\end{align*} 
\end{prop}
\noindent {\bf Proof:}  Assume $ x \in \mathrm{Affine}_{z} $.  For $ y = x $ and $ y = x^*_{z} \; , $  we have the equality (\ref{eqn.bg}), that is,  
\[      c \cdot  \pi(y)  = c \cdot e + \smfrac{1}{1 - \lambda_{\min}(y)} ( z - c \cdot e ) \; . 
\]
Thus,
\begin{align*}
    \frac{ c \cdot  \pi(x) - z^* }{c \cdot e - z^* } & =   \frac{ c \cdot  \pi(x)  - c \cdot  \pi(x^*_{z} )  }{c \cdot e - c \cdot  \pi(x^*_{z} )  } \\
    & = \frac{ \smfrac{1}{1 - \lambdamin(x)} -  \smfrac{1}{1 - \lambdamin(x^*_{z} )} }{ - \smfrac{1}{1 - \lambdamin(x^*_{z} )}} \\  
    & = \frac{ \lambdamin(x^*_{z} ) - \lambdamin(x) }{ 1 - \lambdamin(x)} \; . 
    \end{align*}
Hence,
\begin{gather*}     
 \frac{c \cdot  \pi(x)  - z^* }{c \cdot e - z^* } \leq \epsilon \\
\Leftrightarrow \\ 
   \lambdamin(x^*_{z} ) - \lambdamin(x) \leq \epsilon \, ( 1 - \lambdamin(x)) \\
\Leftrightarrow \\
(1 - \epsilon ) ( \lambdamin(x^*_{z} ) - \lambdamin(x) ) \leq \epsilon ( 1 - \lambdamin(x^*_{z} )) \\
\Leftrightarrow \\ \lambdamin(x^*_{z} ) - \lambdamin(x) \leq \smfrac{\epsilon }{1 - \epsilon } ( 1 - \lambdamin(x^*_{z} )) \; .  
   \end{gather*}
Using Lemma~\ref{lem.bd}    to substitute for the rightmost occurrence of $ \lambda_{\min}(x^*_{z} ) $ completes the proof. \hfill $ \Box $
 \vspace{3mm}
 
In concluding the section, we remark that the basic theory holds for convex conic optimization problems generally.  For example, consider a conic program
\begin{equation}  \label{eqn.bj}
  \left. \begin{array}{rl}
  \min_{x \in {\mathcal E}}  & c \cdot x \\
   \textrm{s.t.} & x \in \mathrm{Affine} \\
    & Ax + b \in {\mathcal K}' \end{array}  \right\} \cp'  
    \end{equation}  
Here, $ A $ is a linear operator from $ {\mathcal E} $ to a Euclidean space $ {\mathcal E}' $, $ b \in {\mathcal E}' $ and $ {\mathcal K}' $ is a proper, closed, convex cone in $ {\mathcal E}' $ with nonempty interior. \vspace{1mm}

For a problem with multiple conic constraints, simply let $ {\mathcal K}' $ be the Cartesian product of the cones. \vspace{1mm}

 Obviously, the optimization problem CP corresponds to the case that $ A $ is the identity, $ b = 0 $ and $ {\mathcal K}' = {\mathcal K}  $.  (Thus, on the surface, $ \mathrm{CP}' $ appears to be more general than CP.) \vspace{1mm}
 
Fix  a feasible point $ e $  for which $ e' := Ae + b \in \int({\mathcal K}') $. Using $ e' $ as the distinguished direction results in a function $ x' \mapsto \lambda_{\min}'(x') $ and a seminorm $ \| \, \, \|_{\infty}' $ on $ {\mathcal E'} $. A seminorm is induced on $ {\mathcal E} $, according to $ v \mapsto \| Av \|_{\infty}' \; , $ for which the closed unit ball centered at $ e $ is the largest subset of $ \{ x: Ax + b \in {\mathcal K}' \} $ with symmetry point $ e $.  The map $ x \mapsto \lambda_{\min}'(Ax + b) $  is Lipschitz continuous:  
\[    | \lambda_{\min}'(Ax + b) - \lambda_{\min}'(Ay+b) | \leq \| A(x-y) \|_{\infty}' \quad \textrm{for all $ x,y \in {\mathcal E} $} \; . \]

If $ x \in \mathrm{Affine} $ satisfies $ c \cdot x < c \cdot e \; , $ then $ \lambda_{\min}'(Ax + b) < 1 $, and the projection of $ x $ (from $ e $) to the boundary of the feasible region is given by 
\[  \pi(x) = e + \smfrac{1}{1 - \lambda_{\min}'(Ax + b)} \, (x - e)  \; . \]
Assuming $ z $ is a scalar satisfying  $ z < c \cdot e $, the problem
\begin{equation}  \label{eqn.bk}
   \begin{array}{rl}
       \max & \lambda_{\min}'(Ax + b) \\
       \textrm{s.t.} & x \in \mathrm{Affine}_{z} \quad \textrm{($ := \{ x \in \mathrm{Affine}: c \cdot x = z \} $)}     \end{array} 
       \end{equation} 
is equivalent to $ \mathrm{CP}' $ in that when $ x \in \mathrm{Affine}_{z} \; , $  $ x $ is optimal for (\ref{eqn.bk})  if and only if $ \pi(x) $ is optimal for $ \mathrm{CP}' $.  Moreover, letting $ x_{z}^* $ denote any optimal solution of (\ref{eqn.bj}), there holds the relation for all $ x \in \mathrm{Affine}_{z} $ and $ 0 < \epsilon < 1 $,  
\[    \frac{c \cdot x - z^* }{ c \cdot e - z^* } \, \leq \, \epsilon \quad \Leftrightarrow \quad \lambda_{\min}'(Ax_{z}^* + b) - \lambda_{\min}'(Ax + b) \, \leq \, \frac{\epsilon}{1 - \epsilon} \, \, \frac{ c \cdot e - z}{ \, \, \,  c \cdot e - z^* } \; , \]
where $ z^* $ is the optimal value of $ \mathrm{CP}' $. \vspace{1mm}

These claims regarding $ \mathrm{CP}' $ are justified with proofs that are essentially identical to the proofs for $ \mathrm{CP} $.  Alternatively, they can be deduced from the results for $ \mathrm{CP} $ by introducing a new variable $ t $ into $ \mathrm{CP}' $ and an equation $ t = 1 $, then replacing $ {\mathcal K}' $ by $ {\mathcal K}  := \{ (x,t): Ax + tb \in {\mathcal K}' \} $, thereby recasting $ \mathrm{CP}' $ to be of the same form as $ \mathrm{CP} $. (Only on the surface does $ \mathrm{CP}' $ appear to be more general than CP.) \vspace{1mm}

We focus on $ \mathrm{CP} $ because notationally its form is least cumbersome. For every result derived in the following sections, an essentially identical result holds for any conic program, even identical in the specific constants. 

\section{ {\bf Applying Supgradient Methods}} \label{sect.c}

In this section we show how the basic theory from Section~\ref{sect.b}  leads to complexity results regarding the solution of the conic program
\[ 
  \left. \begin{array}{rl}
\min  & c \cdot x  \\
\textrm{s.t.} & x \in \mathrm{Affine}  \\
  & x \in {\mathcal K}   \end{array} \right\} \cp 
\] 
  Continue to assume CP has an optimal solution, denote the optimal value by $ z^* $, and let $ e$ be a strictly feasible point, the distinguished direction. 
\vspace{1mm}

Given $ \epsilon > 0 $ and a value satisfying $ z < c \cdot e $, the approach is to apply supgradient methods to approximately solve
\begin{equation}  \label{eqn.ca}
 \begin{array}{rl}
 \max & \lambda_{ \min}(x) \\
\textrm{s.t.} & \mathrm{Affine}_{z}   \; ,  \end{array} 
      \end{equation}
where by ``approximately solve'' we mean that  $ x \in \mathrm{Affine}_{z}  $ is computed for which 
\[ 
\lambda_{\min}(x^*_{z}) - \lambda_{ \min}(x) \,  \leq \,  \frac{\epsilon}{1 - \epsilon } \, \, \frac{c \cdot e - z   }{\, \, \, c \cdot e - z^* }  \; . 
\] 
Indeed, according to Proposition~\ref{prop.be}, the projection $ \pi = \pi(x) $ will then satisfy
\[ 
\frac{c \cdot \pi - z^* }{c \cdot e - z^* } \leq \epsilon  \; . 
\] 

Recall that $ {\mathcal L}  $ is the subspace of which the the affine space $ \mathrm{Affine}_z $ is a  translate.   When supgradient methods are applied to solving (\ref{eqn.ca}), $ {\mathcal L} $  is the subpace onto which supgradients are orthogonally projected. 
\vspace{1mm}

Supgradients and their orthogonal projections depend on the inner product.\footnote{Of course orthogonality in a Euclidean space $ {\mathcal E} $ depends on the chosen inner product ($ \lin u, v \rin = 0 $), and hence so do orthogonal projections.   Recall that supgradients also depend on the chosen inner product, in that for a concave function $ f: {\mathcal E} \rightarrow \mathbb{R} $, the supgradients of $ f $ at $ x  $ are the vectors $ \grad f(x) \in {\mathcal E} $ satisfying
$    f(x) + \lin \grad f(x), v \rin \, \geq \, f(x + v) $ for all $ v \in {\mathcal E} $.}  We allow the ``computational inner product'' to differ from the one relied upon in expressing CP, the inner product written $ u \cdot v $.  We denote the computational inner product by $ \langle \; , \; \rangle $, and let $ \| \, \, \| $ be its norm. \vspace{1mm}

It is an instructive exercise to show that the supdifferential (set of supgradients) at $ x $ for the function $ x \mapsto \lambda_{\min}(x) $ is
\[    \big\{ g : \lin  g, e \rin = 1  \textrm{ and } \big\langle g, y - \big( x - \lambda_{\min}(x) \,  e \big) \big\rangle \geq 0 \textrm{ for all $ y \in {\mathcal K}  $} \} \; , \]
that is, the supdifferential consists of vectors $ g $ such that $ \lin g, e \rin = 1 $ and $ -g $ is in the normal cone to $ {\mathcal K}  $ at $ x - \lambda_{\min}(x) \,  e $. (To begin, note it may be assumed that $ \lambda_{\min}(x) = 0 $, due to (\ref{eqn.ba}).) \vspace{3mm}

For $ z \in \mathbb{R} $, let
\[  M_z := \sup \left\{ \smfrac{| \lambda_{\min}(x) - \lambda_{\min}(y)|}{ \| x - y \| }: x,y \in \mathrm{Affine}_z  \textrm{ and } x \neq y \right\} \; , \]
the Lipschitz constant for the map $ x \mapsto \lambda_{\min}(x) $ restricted to $ \mathrm{Affine}_z $.  Proposition~\ref{prop.ba}  implies $ M_z $ is well-defined (finite), although unlike the Lipschitz constant for the norm appearing there (i.e., $ \| \, \, \|_{\infty} $), $ M_z $ might exceed 1, depending on $ \| \, \, \| $. \vspace{1mm}

We claim the values $ M_z $ are identical for all $ z $.  To see why, consider that for $ z_1 < z_2 < c \cdot e \; , $  a bijection from $ \mathrm{Affine}_{z_1} $ onto $ \mathrm{Affine}_{z_2} $ is provided by the map
\[  x \mapsto y(x) := x + \smfrac{z_2 - z_1}{c \cdot e -  z_1} ( e - x) \; . \]
Observe, using (\ref{eqn.ba}),  
\[  \lambda_{\min}(y(x)) = \smfrac{c \cdot e - z_2}{c \cdot e - z_1} \lambda_{\min}(x) + \smfrac{z_2 - z_1}{c \cdot e - z_1} \; , \]
and thus 
\[  \lambda_{\min}(y(x)) - \lambda_{\min}(y( \bar{x})) = \smfrac{c \cdot e - z_2}{c \cdot e - z_1} \left(   \lambda_{\min}(x) - \lambda_{\min}( \bar{x}) \right) \quad \textrm{for $ x, \bar{x} \in \mathrm{Affine}_{z_1} $} \; . \]
Since, additionally, $ \| y(x) - y( \bar{x}) \| = \smfrac{c \cdot e - z_2}{c \cdot e - z_1} \| x - \bar{x} \| $, it is immediate that the values $ M_z $ are identical for all $ z < c \cdot e $.  A simple continuity argument then implies this value is equal to $ M_{c \cdot e} $. Analogous reasoning shows $ M_z = M_{ c \cdot e} $ for all $ z > c \cdot e $.  In all, $ M_z $ is independent of $ z $, as claimed.  \vspace{1mm}

Let $ M $ denote the common value, i.e., $ M = M_z $ for all $ z $. \vspace{1mm}

The following proposition can be useful in modeling and in choosing the computational inner product (the inner product for whose norm $ M $ is the Lipschitz constant). \vspace{1mm}

Let $ B(e,r) := \{ x: \| x - e \| \leq r \} $.  

\begin{prop} \label{prop.ca} $ M \leq 1/ \bar{r} \; , $ where  $ \bar{r} :=  \max \{ r: B(e,r) \cap \mathrm{Affine}_{c \cdot e} \subseteq {\mathcal K}  \} $ 
\end{prop}
\noindent {\bf Proof:} According to Proposition~\ref{prop.ba}, 
\[   | \lambda_{\min}(x) - \lambda_{\min}(y)| \leq \| x - y \|_{\infty} \quad \textrm{for all $ x, y $} \; . \]
Consequently, it suffices to show $ \| x - y \| \leq \bar{r} \| x - y \|_{\infty} $ for all $ x,y \in \mathrm{Affine}_{c \cdot e} \; , $ i.e., it suffices to show for all $  v \in {\mathcal L} $ that $ \| v \| \leq \bar{r} \| v \|_{\infty} $. \vspace{1mm}

 However, according to the discussion just prior to Proposition~\ref{prop.ba},  $ B_{\infty}(e,1) $ is the largest set which both is contained in $ {\mathcal K} $ and has symmetry point $ e $, from which follows that $ B_{\infty}(e,1) \cap \mathrm{Affine}_{c \cdot e} $ is the largest set which is both contained in $ {\mathcal K}  \cap \mathrm{Affine}_{c \cdot e} $ and has symmetry point $ e $.  Hence 
 \[  B(e, \bar{r}) \cap \mathrm{Affine}_{ c \cdot e} \subseteq B_{\infty}(e,1) \cap \mathrm{Affine}_{c \cdot e} \; , \]
implying $ \| v \| \leq \bar{r} \| v \|_{\infty} $ for all $ v \in {\mathcal L} $. \hfill $ \Box $ \vspace{3mm}

In passing we remark that in the context of $ \mathrm{CP}' $ -- the conic program (\ref{eqn.bj})  -- the role of $ \bar{r}$ in the above proposition is played by  
\[ \bar{r} = \max \left\{ r: B(e, r) \cap \mathrm{Affine}_{ c \cdot e} \subseteq \{ x: Ax + b \in {\mathcal K}' \} \right\}\; . \]
Then, $ | \lambda_{\min}(Ax + b) - \lambda_{\min}(Ay+b)| \leq (1/ \bar{r}) \| x - y \| $ for all $ x,y \in \mathrm{Affine}_z $ and $ z \in \mathbb{R}$.   \vspace{1mm}

Towards considering specific supgradient methods, we recall the following standard and elementary result, rephrased for our setting:

\begin{lemma}   \label{lem.cb} 
  Assume $ z \in \mathbb{R} $ and $ x, y \in \mathrm{Affine}_z \; .  $ Let $ g $ be the projection of a supgradient $ \grad \lambda_{\min}(x) $ onto $ {\mathcal L} $ (the subspace of which $ \mathrm{Affine}_z $ is a translate).  

For all scalars $ \alpha  $,
\begin{equation}  \label{eqn.cb} 
   \| ( x + \alpha g) - y \|^2 \leq \| x - y \|^2 - 2 \alpha  \left(   \lambda_{\min}(y) - \lambda_{\min}(x) \right) + \alpha^2 \| g \|^2 \; . \end{equation}
\end{lemma} 
\noindent {\bf Proof:}  Simply observe
\begin{align*}
\| ( x + \alpha g) - y \|^2 & = \| x - y \|^2 + 2 \alpha \lin g, x - y \rin + \alpha^2 \| g \|^2 \\
& =  \| x - y \|^2 - 2 \alpha \lin \grad \lambda_{\min}(x), y - x \rin + \alpha^2 \| g \|^2 \quad \textrm{(by $ x - y \in {\mathcal L} $)} \\
& \leq \| x - y \|^2 - 2 \alpha \left(   \lambda_{\min}(y) - \lambda_{\min}(x) \right)  + \alpha^2 \| g \|^2 \; , 
\end{align*}
the inequality due to concavity of the map $ x \mapsto \lambda_{\min}(x) \; . $ \hfill $ \Box $ \vspace{3mm}

We present and analyze two algorithms. We begin by considering the ideal case in which $ z^* $, the optimal value for CP, is known.  The main result for the algorithm here provides a benchmark to which to compare our result for the general algorithm  developed subsequently. \vspace{1mm}

Knowing $ z^* $ is not an entirely implausible situation.  For example, if strict feasibility holds for a primal conic program and for its dual
\[ \begin{array}{rl}
          \min & \bar{c}^T x \\
          \textrm{s.t.} & \bar{A}x = \bar{b} \\
     & x \in \bar{{\mathcal K} } \end{array} \qquad  \begin{array}{rl}
          \max & \bar{b}^T y \\
           \textrm{s.t.} & \bar{A}^T y + s = \bar{c} \\
                    & s \in \bar{{\mathcal K} }^* \; , \end{array} \]
then the combined primal-dual conic program is known to have optimal value equal to zero:
\[   \begin{array}{rl}
       \min & \bar{c}^T x - \bar{b}^T y \\
       \textrm{s.t.} & \bar{A} x = \bar{b} \\
                    & \bar{A}^T y + s = \bar{c} \\
                    & (x,s) \in \bar{{\mathcal K} } \times \bar{{\mathcal K} }^* \; . \end{array} \]

 \noindent 
\hrulefill
                    
\noindent {\bf Algorithm 1:} \vspace{1mm} 

\noindent  
(0) Input: $ z^* $, the optimal value of CP, \\
$ \textrm{~} $ \qquad \qquad  \,  $ e $, a strictly feasible point for CP, and \\
$ \textrm{~} $ \qquad \qquad  \,    $ \bar{x}  \in \mathrm{Affine} $ satisfying $ c \cdot \bar{x}  < c \cdot e \; . $  \\
$ \textrm{~} $ \quad Initialize: Let $ x_0 = e + \frac{c \cdot e - z^*}{c \cdot e - c \cdot \bar{x} }(\bar{x} -e) $ \,  (thus, $ c \cdot x_0 = z^* $), \\
$ \textrm{~} $ \qquad \qquad  \qquad  \quad    and let $ \pi_0 = \pi(x_0) $ ($ = \pi(\bar{x} ) $). \\
(1) Iterate: $ x_{k+1} = x_k - \smfrac{\lambda_{\min}(x_k)}{\| g_k \|^2} g_k $, \\
$ \textrm{~} $ \qquad \qquad \qquad  \qquad   where $ g_k $ is the projection of a supgradient $ \grad \lambda_{\min}(x_k) $ onto $ {\mathcal L} $. \\
$ \textrm{~} $ \qquad \qquad  \qquad  Let $ \pi_{k+1} = \pi(x_{k+1}) $ \, \, (which is feasible for CP)  \vspace{-1.5mm}

\noindent 
\hrulefill
\vspace{2mm}

All of the iterates $ x_k $ lie in $ \mathrm{Affine}_{z^*} $, and hence, $ \lambda_{\min}(x_k) \leq 0 $, with equality if and only if $ x_k $ is feasible (and optimal) for CP. \vspace{1mm}

For all scalars $ z < c \cdot e $ and for $ x \in \mathrm{Affine}_z \; ,  $ define
\[  
   \mathrm{dist}_z(x) :=  \min \{ \| x - x_z^* \| : x_z^* \textrm{ is optimal for  (\ref{eqn.ca})} \} \; . \] 

\begin{prop}  \label{prop.cc} 
The iterates for Algorithm 1 satisfy
\[ 
 \max \{ \lambda_{\min}(x_k) : k = \ell, \ldots, \ell + m \} \geq -M \, \mathrm{dist}_{z^*}(x_{\ell})/ \sqrt{m+1} \; . 
\] 
\end{prop} 
\noindent {\bf Proof:}  Lemma~\ref{lem.cb}  implies
\begin{align*} 
   & \mathrm{dist}_{z^*}(x_{k+1})^2 \\ & \leq  \mathrm{dist}_{z^*}(x_k)^2 - 2 (-\lambda_{\min}(x_k)/\| g_k \|^2) \,  ( 0 - \lambda_{\min}(x_k)) \, + \, (\lambda_{\min}(x_k)/\| g_k \|)^2 \\  
    & = \mathrm{dist}_{z^*}(x_k)^2 - (\lambda_{\min}(x_k)/\| g_k \|)^2 \; , 
    \end{align*}
and thus by induction (and using $ \| g_k \| \leq M $),
\begin{align*}
\mathrm{dist}_{z^*}(x_{\ell + m + 1 })^2 & \leq \mathrm{dist}_{z^*}(x_{\ell})^2 - \sum_{k= \ell }^{\ell + m} ( \lambda_{\min}(x_k)/M)^2  \\
& \leq  \mathrm{dist}_{z^*}(x_{\ell})^2 - \smfrac{m + 1}{M^2}  \min \{ \lambda_{\min}(x_k)^2 : k = \ell, \ldots, \ell + m \} \; , 
\end{align*}
implying the proposition (keeping in mind $ \lambda_{\min}(x_k) \leq 0 $). \hfill $ \Box $
 \vspace{3mm}
  
We briefly digress to consider the case of $ {\mathcal K}  $ being polyhedral, where already an interesting result is easily proven. The following corollary is offered only as a curiosity, as the constant $ C_2 $  typically is so large as to render the lower bound on $ \ell$ meaningless except for minuscule $ \epsilon$.
 \vspace{1mm}

\begin{cor}  \label{cor.cd} 
Assume $ {\mathcal K}  $ is polyhedral. There exist constants $ C_1 $ and $ C_2 $ (dependent on CP, $ e $, $ \bar{x} $ and the computational inner product), such that for all $ 0 < \epsilon < 1 $, 
\[  \ell  \geq C_1 + C_2 \log(1/ \epsilon ) \quad \Rightarrow \quad \min_{k \leq \ell} \frac{c \cdot \pi_k - z^*}{c \cdot e - z^*} \,  \leq \, \epsilon \; . \]
\end{cor}
\vspace{2mm}

  For first-order methods, such a logarithmic bound in $ \epsilon $  was initially established by Gilpin, Pe\~{n}a and Sandholm \cite{gilpin2012first}. They did not assume an initial feasible point $ e $ was known, but neither did they require the computed solution to be feasible (instead, constraint residuals were required to be small).  They relied on an accelerated gradient method, along with the smoothing technique of Nesterov \cite{nesterov2005smooth}.  As is the case for the above result, they assumed the optimal value of CP to be known apriori, and they restricted $ {\mathcal K}  $ to be polyhedral. \vspace{2mm}
  
The proof of the corollary depends on the following simple lemma. \vspace{1mm}        

\begin{lemma}  \label{lem.ce} 
For Algorithm 1, the iterates satisfy
\[  \frac{c \cdot \pi_k - z^*}{c \cdot e - z^*} = \frac{- \lambda_{\min}(x_k) }{1 - \lambda_{\min}(x_k)} \; . \]
\end{lemma}
\noindent {\bf Proof:}  Immediate from $ \pi_k = e + \frac{1}{1 - \lambda_{\min}(x_k)} (x_k - e) $ and $ c \cdot x_k = z^* $. \hfill $ \Box $
 \vspace{3mm}  
 
\noindent {\bf Proof of Corollary~\ref{cor.cd}:} 
 With $ {\mathcal K}  $ being polyhedral, the concave function $ x \mapsto \lambda_{\min}(x) $ is piecewise linear, and thus there exists a positive constant $ C $ such that
\[      \mathrm{dist}_{z^*}(x) \leq - C \, \lambda_{\min}(x) \quad \textrm{for all $ x \in \mathrm{Affine}_{z^*} $} \; . \]
Then Proposition~\ref{prop.cc}  gives
\[  
 \max  \{ \lambda_{\min}(x_k) : k = \ell, \ldots, \ell + m \} \geq C \, M \, \lambda_{\min}(x_{\ell})/ \sqrt{m+1} \; , \]
from which follows 
 \[  \max  \{ \lambda_{\min}(x_k) : k = \ell, \ldots, \ell + \lceil (2CM)^2 \rceil \} \geq \smfrac{1}{2} \lambda_{\min}(x_{\ell}) \; , \]
i.e., $ \lambda_{\min}(x_{\ell}) $ is ``halved'' within $ \lceil (2CM)^2 \rceil $ iterations.  The proof is easily completed using Lemma~\ref{lem.ce}.   \hfill $ \Box $ \vspace{3mm}

We now return to considering general convex cones $ {\mathcal K}  $. \vspace{1mm}

The iteration bound provided by Proposition~\ref{prop.cc}  bears little obvious connection to the geometry of the conic program CP, except in that the constant $ M $ is related to the geometry by Proposition~\ref{prop.ca}.  The other constant -- $ \mathrm{dist}_{z^*}(x_k) $ -- does not at present have such a clear geometrical connection to CP.  We next observe a meaningful connection.  \vspace{1mm}

 The {\em level sets}  for CP are the sets
\[  \mathrm{Level}_{z} =   \mathrm{Affine}_{z} \cap {\mathcal K}   \; ,   \]
that is, the largest feasible sets for CP on which the objective function is constant\footnote{There is possibility of confusion here, as in the optimization literature, the terminology ``level set'' is often used for the portion of the feasible region on which the (convex) objective function does not exceed a specified value rather than -- as for us -- exactly equals the value.  Our choice of terminology is consistent with the general mathematical literature, where the region on which a function does not exceed a specified value is referred to as a sublevel set, not a level set.}.  If $ z < z^* $, then $ \mathrm{Level}_{z} = \emptyset \; . $    \vspace{1mm}

If some level set is unbounded, then either CP has unbounded optimal value or can be made to have unbounded value with an arbitrarily small perturbation of $ c $.  Thus, in developing numerical optimization methods, it is natural to focus on the case that level sets for CP are bounded.  \vspace{1mm}  
     
For scalars $ z  $, let 
\[    \diamval  :=  \sup  \{  \| x - y \|: x,y \in \mathrm{Level}_{z}   \} \; ,    \]
the diameter of $ \mathrm{Level}_{z} $.   (If $ \mathrm{Level}_z = \emptyset $, let $ \diamval := - \infty  $.) \vspace{1mm}

\begin{lemma}  \label{lem.cf} 
Assume $ x \in \mathrm{Affine}_{z^*} $, and let $ \pi = \pi(x) $. Then
       \[  \mathrm{dist}_{z^*}(x) \,  = \,  (1 - \lambda_{\min}(x)) \, \mathrm{dist}_{c \cdot \pi}(\pi) \,  = \, \frac{ \mathrm{dist}_{c \cdot \pi}(\pi)}{1 - \frac{c \cdot \pi - z^*}{c \cdot e - z^*}} \,  \leq \, \frac{ \mathrm{diam}_{c \cdot \pi}}{1 - \frac{c \cdot \pi - z^*}{c \cdot e - z^*}} \; .  \] 
 \end{lemma}
\noindent {\bf Proof:}  Since 
\begin{equation}  \label{eqn.cc} 
  \pi = e + \smfrac{1}{1 - \lambda_{\min}(x)}(x-e) \; , \end{equation}
  Theorem~\ref{thm.bc}  implies that the maximizers of the map $ y \mapsto \lambda_{\min}(y) $ over $ \mathrm{Affine}_{ c \cdot \pi} $ are precisely the points of the form
\[   x_{c \cdot \pi}^* =   e + \smfrac{1}{1 - \lambda_{\min}(x)}(x_{z^*}^*-e) \; , \]
where $ x_{z^*}^* $ is a maximizer of the map when restricted to $ \mathrm{Affine}_{z^*} $ (i.e., is an optimal solution of CP).  Observing  
\[    \pi - x_{c \cdot \pi}^* = \smfrac{1}{1 - \lambda_{\min}(x)}  (x - x_{z^*}^*) \; , \]
it follows that 
\[   \mathrm{dist}_{c \cdot \pi}( \pi) =  \smfrac{1}{1 - \lambda_{\min}(x)}  \, \mathrm{dist}_{z^*}(x) \; , \]
establishing the first equality in the statement of the lemma.  The second equality then follows easily from (\ref{eqn.cc}) and $ c \cdot x = z^* $. The inequality is due simply to $ \pi, x_{c \cdot \pi}^* \in \mathrm{Level}_{ c \cdot \pi} $, for all optimal solutions $ x_{c \cdot \pi}^* $ of the CP-equivalent problem (\ref{eqn.ca}) (with $ z = c \cdot \pi $).   \hfill $ \Box $
 \vspace{3mm}
 
For scalars $ z $, define
\[  \mathrm{Diam}_z := \max \{ \mathrm{diam}_{z'}: z' \leq z \} \; , \]
the ``horizontal diameter'' of the sublevel set consisting of points $ x $ that are feasible for CP and satisfy $ c \cdot x \leq z $.  For $ z^* < z < c \cdot e $, the  value $ \mathrm{Diam}_{z} $ can be thought of as a kind of condition number for CP, because $  \mathrm{Diam}_{z} $ being large is an indication that the optimal value for CP is relatively sensitive to perturbations in the objective vector $ c $.  (Related quantities have played the role of condition number in the complexity theory of interior-point methods  -- c.f., \cite{freund2004complexity}.) \vspace{1mm}

For $ z^* \leq  z < c \cdot e $, define
\[     \mathrm{Dist}_z := \sup \{ \mathrm{diam}_{z'}(x): z' \leq z \textrm{ and } x \in \mathrm{Level}_{z'} \}  \; . \]
Clearly, there holds the relation
\[    \mathrm{Dist}_z \leq \mathrm{Diam}_z \; , \]
and hence if the ``condition number'' $ \mathrm{Diam}_z $ is only of modest size, so is the value $ \mathrm{Dist}_z $.  \vspace{1mm}

Following is our main result for Algorithm 1. By substituting $  \mathrm{Diam}_{c \cdot \pi_0} $ for $  \mathrm{Dist}_{ c \cdot \pi_0} $, and $ 1/\bar{r}$ for $ M $ (where $ \bar{r} $ is as in Proposition~\ref{prop.ca}), the statement of the theorem becomes phrased in terms clearly reflecting the geometry of CP. \vspace{1mm}

\begin{thm}  \label{thm.cg} 
Assume $ 0 < \epsilon < \frac{c \cdot \pi_0 - z^*}{c \cdot e - z^*} $, where $ \pi_0 = \pi(x_0) $ is the initial $ \mathrm{CP} $-feasible point for Algorithm 1 (i.e., assume $ \pi_0 $ does not itself satisfy the desired accuracy).  Then 
\begin{align*}
 &   \ell  \geq (2M \, \mathrm{Dist}_{c \cdot \pi_0})^2 \, \, \left( \, \frac{4}{ 3  } \left( \frac{1 - \epsilon }{ \epsilon } \right)^2 + 4 \left(   \frac{1 - \epsilon}{\epsilon} \right)  \right. \\
& \qquad \qquad \qquad  \qquad  \qquad  \qquad   \left. + \log_2 \left( \frac{\frac{c \cdot \pi_0 - z^*}{c \cdot e - z^*} }{\epsilon} \right) + \log_2 \left( \frac{1 - \epsilon }{1 - \frac{c \cdot \pi_0 - z^*}{c \cdot e - z^*} } \right) \, + 1 \, \right)       \\ &  \qquad        \quad \Rightarrow \quad \min_{k \leq \ell } \,   \frac{c \cdot \pi_k - z^*}{c \cdot e - z^*} \, \leq \, \epsilon \; . \end{align*} 
\end{thm}
\noindent {\bf Proof:}   To ease notation, let $  \lambda_k := \lambda_{\min}(x_k) \; . $ \vspace{1mm}

 Let $ k_0 = 0 $ and recursively define $ k_{i+1} $ to be the first index for which $ \lambda_{k_{i+1}} \leq \lambda_{k_i}/2 $ (keeping in mind $ \lambda_k \leq 0 $ for all $ k $).  Proposition~\ref{prop.cc} implies  
\begin{align}
   k_{i+1} - k_i + 1   & \leq \left( \frac{2M \, \mathrm{dist}_{z^*}(x_{k_i})}{\lambda_{k_i}} \right)^2 \nonumber \\
   & \leq \left( \, 2 M \, \mathrm{dist}_{c \cdot \pi_{k_i}}(\pi_{k_i}) \, \frac{ 1 - \lambda_{k_i}}{ \lambda_{k_i}} \,  \right)^2 \quad \textrm{(by Lemma~\ref{lem.cf})} \nonumber \\ 
  & \leq   \left( \, 2 M \, \mathrm{Dist}_{c \cdot \pi_0} \, \frac{ 1 - \lambda_{k_i}}{ \lambda_{k_i}} \,  \right)^2 \; ,  \label{eqn.cd} 
  \end{align}
  where the final inequality is due to $ c \cdot \pi_{k_i} $ ($ i = 0, 1, \ldots $) being a decreasing sequence (using Lemma~\ref{lem.ce}). \vspace{1mm}
  
Let $ i' $ be the first sub-index for which $ \lambda_{k_{i'}} \geq - \epsilon/(1 - \epsilon) $. Lemma~\ref{lem.ce}  implies
\[   \frac{c \cdot \pi_{k_{i'}} - z^*}{ c \cdot e - z^* } \, \leq  \, \epsilon \; , \]
Thus, to prove the theorem, it suffices to show $ \ell  = k_{i'} $ satisfies the inequality in the statement of the theorem.  \vspace{1mm}

Note $ i' > 0 $ (because, by assumption, $ \epsilon < \frac{c \cdot \pi_0 - z^*}{c \cdot e - z^*} $). Observe, then, 
\begin{align}
    i' &  < 1 +   \log_2 \left( \frac{ \lambda_0}{- \epsilon/(1 - \epsilon )} \right)  \nonumber \\
       &  = 1 +   \log_2 \left( \frac{\frac{c \cdot \pi_0 - z^*}{c \cdot e - z^*} }{ \epsilon} \right) + \log_2 \left( \frac{1 - \epsilon }{1 - \frac{c \cdot \pi_0 - z^*}{c \cdot e - z^*} } \right)    \label{eqn.ce} 
       \end{align}  
       (again using Lemma~\ref{lem.ce}).     
  \vspace{1mm}

Additionally,
\begin{align*}
k_{i'} & = \sum_{i=0}^{i'-1} k_{i+1} - k_i  \\ &
    \leq (2M \, \mathrm{Dist}_{c \cdot \pi_0})^2 \, \sum_{i=0}^{i'-1} \left( \frac{1 - \lambda_{k_i}}{\lambda_{k_i}} \right)^2 \quad \textrm{(by (\ref{eqn.cd})} \\
       & \leq (2M \, \mathrm{Dist}_{c \cdot \pi_0})^2 \,  \sum_{i=0}^{i'-1} \left( \frac{1 - 2^i \lambda_{k_{i'-1}}}{2^i \lambda_{k_{i'-1}}} \right)^2  \\
       & \leq  (2M \, \mathrm{Dist}_{c \cdot \pi_0})^2 \, \sum_{i=0}^{i'-1} \left( \frac{1 + 2^i \epsilon /(1 - \epsilon)}{2^i \epsilon /(1 - \epsilon) } \right)^2 \\
     &  = (2M \, \mathrm{Dist}_{c \cdot \pi_0})^2 \, \sum_{i=0}^{i'-1} \left( 1 + \frac{1}{2^i} \, \frac{1 - \epsilon}{\epsilon } \right)^2 \\    
& \leq (2M \, \mathrm{Dist}_{c \cdot \pi_0})^2 \, \left(  i' + 4  \frac{1- \epsilon }{ \epsilon }  + \frac{4}{3} \left( \frac{1- \epsilon }{ \epsilon } \right)^2  \right) \; .    
\end{align*}
Using (\ref{eqn.ce})  to substitute for $ i' $ completes the proof. \hfill $ \Box $
\vspace{3mm}

For our second algorithm, we discard the requirement of knowing $ z^* $, the optimal value for CP.  Now we require that $ \epsilon $ (the desired relative-accuracy) be input. \vspace{2mm} 

\noindent 
\hrulefill
 
\noindent {\bf Algorithm 2:} \vspace{1mm}

\noindent   
(0) Input: $ 0 < \epsilon < 1 $ ,\\
 $ \textrm{~} $ \qquad \qquad  \,  $ e $, a strictly feasible point for CP, and \\
$ \textrm{~} $ \qquad \qquad  \,   $ \bar{x}  \in \mathrm{Affine} $ satisfying $ c \cdot \bar{x}  < c \cdot e $.  \\
$ \textrm{~} $ \quad Initialize: $ x_0 = \pi_0 = \pi(\bar{x} ) $ \\
(1) Iterate: $ \tilde{x}_{k+1} :=  x_k + \smfrac{ \epsilon }{2 \| g_k \|^2} g_k $, \\
$ \textrm{~} $ \qquad \qquad \qquad  \qquad   where $ g_k $ is the projection of a supgradient $ \grad \lambda_{\min}(x_k) $ onto $ {\mathcal L} $. \\
$ \textrm{~} $ \qquad \qquad  \quad Let $ \pi_{k+1} :=  \pi(\tilde{x}_{k+1}) \; . $ \\
$ \textrm{~} $ \qquad \qquad  \quad If $ c \cdot (e - \pi_{k+1}) \geq \smfrac{4}{3}  \, c \cdot (e - \tilde{x}_{k+1}) $, then let $ x_{k+1} = \pi_{k+1} $; \\
$ \textrm{~} $ \qquad \qquad  \qquad  \qquad  \qquad  \qquad  \qquad  \qquad  \qquad  \qquad  \qquad   else, let $ x_{k+1} = \tilde{x}_{k+1} \; . $   \vspace{-1.5mm}

\noindent 
\hrulefill
\vspace{2mm}

Unsurprisingly, the iteration bound we obtain for Algorithm 2 is worse than our result for Algorithm 1, but perhaps surprisingly, the bound is not excessively worse, in that the factor for $ 1/ \epsilon^2 $ is essentially unchanged (it's the factor for $ 1/\epsilon $ that increases, although typically not by an excessive amount). \vspace{2mm}

\begin{thm}  \label{thm.ch} 
 Assume $ 0 < \epsilon \leq \frac{c \cdot \pi_0 - z^*}{c \cdot e - z^*}  $.  For the iterates of Algorithm 2,
\begin{align*}
 &   \ell  \geq 8 \, (M \, \mathrm{Dist}_{c \cdot \pi_0} )^2 \, \left( \, \frac{1}{\epsilon^2} \, + \, \frac{1}{\epsilon} \, \log_{4/3} \left( \frac{ 1 - \epsilon }{1 - \frac{c \cdot \pi_0 - z^*}{c \cdot e - z^*}} \right) \, + \, 1 \,  \right)              \\ &  \qquad        \quad \Rightarrow \quad \min_{k \leq \ell} \,   \frac{c \cdot \pi_k - z^*}{c \cdot e - z^*} \, \leq \, \epsilon \; . \end{align*} 
\end{thm}
\noindent {\bf Proof:} In order to distinguish the iterates obtained by projecting to the boundary, we record a notationally-embellished rendition of the algorithm which introduces a distinction between ``inner iterations'' and ``outer iterations'':
\vspace{2mm}

\noindent 
\hrulefill
 
\noindent  {\em Algorithm 2  (notationally-embellished version):} \vspace{1mm}

\noindent  
(0) Input: $ 0 < \epsilon < 1 \; , $ $ e $ and $ \bar{x} $.  \\ 
$ \textrm{~} $ \quad Initialize: $ y_{1,0} = \pi(\bar{x} ) \; , $  \\
$ \textrm{~} $ \qquad \qquad \qquad   $ i = 1$ (outer iteration counter), \\
$ \textrm{~} $ \qquad \qquad  \qquad     $ j = 0 $ (inner iteration counter). \\
(1) Compute $ y_{i,j+1} = y_{i,j} + \smfrac{\epsilon }{2 \| g_{i,j} \|^2} \, g_{i,j} \; , $ \\
$ \textrm{~} $ \qquad \qquad  \qquad   where $ g_{i,j} $ is the projection of a supgradient $ \grad \lambda_{\min}(y_{i,j}) $ onto $ {\mathcal L} $. \\
(2) If $ c \cdot (e - \pi(y_{i,j+1})) \geq \smfrac{4}{3} \, c \cdot (e - y_{i,j+1} ) \; , $ \\
   $ \textrm{~} $ \qquad \qquad  then let \, $ y_{i+1,0} = \pi(y_{i,j+1}) $, \,  $ i + 1 \rightarrow i $ \, and  \, $ 0 \rightarrow j \; $;\\
$ \textrm{~} $ \quad \quad \quad   else, let \, $ j+1 \rightarrow j \; . $  \\
(3) Go to step 1.  \vspace{-1.5mm}

\noindent 
\hrulefill
\vspace{2mm}

For each outer iteration $ i $, all of the iterates $ y_{i,j} $ have the same objective value.  Denote the value by $ z_i $.  Obviously, $ z_1 $ is equal to the value $ c \cdot \pi_0 $ appearing in the statement of the theorem.  Let
\[ \mathrm{Dist} := \mathrm{Dist}_{c \cdot \pi_0} \quad = \mathrm{Dist}_{z_1} \; . \]

Step 2 ensures     
\begin{equation}  \label{eqn.cf} 
    c \cdot e - z_{i+1} \geq \smfrac{4}{3}  ( c \cdot e - z_i)  \; . 
\end{equation}  
Thus, $ z_1, z_2, \ldots $ is a strictly decreasing sequence.  Consequently, as $ y_{i,0} \in \mathrm{Level}_{z_i} $, we have $  \mathrm{dist}_{z_i}(y_{i,0}) \leq \mathrm{Dist} $ for all $ i $.  \vspace{1mm}

From (\ref{eqn.cf})  we find for scalars $ \delta > 0 $ that
\[  
 \frac{c \cdot e - z_{i+1}}{ c \cdot e - z^*} <  \delta \quad \Rightarrow \quad i < \log_{4/3} \left( \frac{\delta}{ \frac{c \cdot e - z_1}{c \cdot e - z^*}} \right) \quad = \log_{4/3} \left( \frac{\delta}{ 1 - \frac{z_1 - z^*}{c \cdot e - z^*}} \right) \; , 
\]
and thus, for $ \epsilon < 1 $, 
\begin{equation}   \label{eqn.cg} 
  \frac{z_i - z^*}{c \cdot e - z^*} >  \epsilon \quad \Rightarrow \quad  i <  1 + \log_{4/3} \left( \frac{1 - \epsilon}{ 1 - \frac{z_1 - z^*}{c \cdot e - z^*}} \right) \; . 
  \end{equation}
Hence, if an outer iteration $ i $ fails to satisfy the inequality on the right, the initial inner iterate $ y_{i,0} $ fulfills the goal of finding a CP-feasible point $ \pi $  satisfying $ \frac{c \cdot \pi - z^*}{c \cdot e - z^*} \leq \epsilon $ (i.e., the algorithm has been successful no later than the start of outer iteration $ i $).  Also observe that (\ref{eqn.cg})  provides (letting $ \epsilon \uparrow 0 $) an upper bound on $ I $, the total number of outer iterations:
\begin{equation}  \label{eqn.ch} 
          I \leq 1 + \log_{4/3} \left( \frac{1}{ 1 - \frac{z_1 - z^*}{c \cdot e - z^*}} \right) \; . 
  \end{equation}
  \vspace{1mm}

For $ i = 1, \ldots, I $,  let $ J_i $ denote the number of inner iterates computed during outer iteration $ i $, that is, $ J_i $ is the largest value $ j $ for which $ y_{i,j} $ is computed. Clearly, $ J_I = \infty $, whereas $ J_1, \ldots, J_{I-1} $ are finite. \vspace{1mm}

To ease notation, let $ \lambda_{i,j} := \lambda_{\min}(y_{i,j}) $, and let $ \lambda_i^* := \lambda_{\min}(x_{z_i}^*) $, the optimal value of
\[ \begin{array}{rl}
    \max & \lambda_{\min}(x) \\
     \textrm{s.t.} & x \in \mathrm{Affine}_{z_i} \; . \end{array} \]
According to Lemma~\ref{lem.bd},
\begin{equation}  \label{eqn.ci} 
\lambda_i^* = \frac{z_i - z^*}{c \cdot e - z^*} \; . 
\end{equation}
It is thus valid, for example, to substitue $ \lambda_i^* $ for $ \frac{z_i - z^*}{c \cdot e - z^*} $ in (\ref{eqn.cg}).   Additionally, (\ref{eqn.ci})  implies (\ref{eqn.cf})  to be equivalent to
\begin{equation}  \label{eqn.cj} 
   1 - \lambda_{i+1}^* \geq \smfrac{4}{3} ( 1 - \lambda_i^* ) \; . 
   \end{equation} 
         
For any point $ y $, we have $ \pi(y) = e + \smfrac{1}{1 - \lambda_{\min}(y)} ( y - e) $,  and thus,
\[   \frac{ c \cdot e - c \cdot \pi(y)}{c \cdot e - c \cdot y} = \frac{1}{1 - \lambda_{\min}(y)} \; . \]
Hence, 
\[ \frac{c \cdot e - c \cdot \pi(y)}{c \cdot e - c \cdot y} \geq \frac{4}{3} \quad \Leftrightarrow \quad \lambda_{\min}(y) \geq 1/4 \; . \]
Consequently,   
\begin{equation}  \label{eqn.ck} 
          \lambda_{i,j} < 1/4 \, \, \textrm{ for $ j < J_i $} \; . 
         \end{equation} 
   
We use the following relation implied by Lemma~\ref{lem.cb}: 
\begin{equation}  \label{eqn.cl} 
 \mathrm{dist}_{z_i}(y_{i,j+1})^2 \leq  \mathrm{dist}_{z_i}(y_{i,j})^2 - \smfrac{\epsilon }{\| g_{i,j} \|^2 } ( \lambda_i^* - \lambda_{i,j} ) + (  \smfrac{\epsilon}{2 \| g_{i,j} \| })^2 \; .
 \end{equation}

We begin bounding the number of inner iterations by showing
 \begin{equation}  \label{eqn.cm} 
       \lambda_i^* \geq \max \{ \smfrac{1}{2} , \epsilon \}  \quad \Rightarrow \quad  J_i \leq \frac{8 (M \, \mathrm{Dist})^2}{ \epsilon } \; .   
       \end{equation}
Indeed, for $ j < J_i \; , $  
\begin{align*}
 &    - \epsilon \,  ( \lambda_i^* - \lambda_{i,j} ) + \smfrac{1}{4} \,  \epsilon^2 \\
& <  - \epsilon \, ( \max \{ \smfrac{1}{2}, \epsilon \}  - \smfrac{1}{4} ) + \smfrac{1}{4} \epsilon^2 \quad \textrm{(using (\ref{eqn.ck})} \\
& = \min  \left\{ \smfrac{1}{4} ( \epsilon^2 - \epsilon), \smfrac{1}{4} \epsilon - \smfrac{3}{4}  \epsilon^2  \right\}   \\
& \leq \smfrac{3}{4} \,   \smfrac{1}{4} ( \epsilon^2 - \epsilon) + \smfrac{1}{4} ( \smfrac{1}{4} \epsilon - \smfrac{3}{4}  \epsilon^2 ) \\
& = - \smfrac{1}{8} \epsilon \; .     
\end{align*} 
Thus, according to (\ref{eqn.cl}), for $ j < J_i $,
\[  \mathrm{dist}_{z_i}(y_{i,j+1})^2 \leq  \mathrm{dist}_{z_i}(y_{i,j})^2 - \frac{ \epsilon}{8  M^2} \; , \]
inductively giving
\begin{align*}
    \mathrm{dist}_{z_i}(y_{i,j+1})^2 & \leq \mathrm{dist}_{z_i}(y_{i,0})^2 - \frac{(j+1) \, \epsilon}{8 \,  M^2 } \\
        & \leq \mathrm{Dist}^2 - \frac{ (j+1) \, \epsilon}{ 8  M^2} \; . 
        \end{align*}
  The implication (\ref{eqn.cm})  immediately follows. \vspace{1mm}      

The theorem is now readily established in the case $ \epsilon \geq 1/2 $. Indeed, because of the identity (\ref{eqn.ci}), the quantity on the right of (\ref{eqn.cg})  provides an upper bound on the number of outer iterations $ i $ for which $ \lambda_i^* > \epsilon $, whereas the quantity on the right of (\ref{eqn.cm})  gives, assuming $ \epsilon \geq 1/2 $, an upper bound on the number of inner iterations for each of these outer iterations.  However, for the first outer iteration satisfying $ \lambda_i^* \leq \epsilon $, the initial iterate $ y_{i,0} $ ($ = \pi(y_{i,0})) $ itself achieves the desired accuracy $ \frac{c \cdot \pi - z^*}{c \cdot e - z^*} \leq \epsilon $. Thus, the total number of inner iterations made before the algorithm is successful is at most the product of the two quantities, which is seen not to exceed the iteration bound in the statement of the theorem. \vspace{1mm}

It remains to consider the case $ \epsilon < 1/2 $.  \vspace{1mm}

For any outer iteration $ i $  for which $ \lambda_i^* < 3/4 $ (and for any $ 0 < \epsilon < 1 $),  let
\[  \widehat{J}_i := \left\lceil  \frac{1}{\smfrac{3}{4}  - \lambda_i^*} \,  \left( \frac{M \, \mathrm{Dist}}{  \epsilon } \right)^2 - 1 \right\rceil  \; . \]
We claim that either
\begin{equation}  \label{eqn.cn} 
    J_i \leq   \widehat{J}_i \quad \textrm{or} \quad  \min \left\{ \frac{c \cdot \pi(y_{i, j }) - z^*}{ c \cdot e - z^* }: j = 0, \ldots, \widehat{J}_i \right\} \, \leq \, \epsilon    \; . 
\end{equation}
Consequently, if $ J_i > \widehat{J}_i $,  the algorithm will achieve the goal of computing a point $ y $ satisfying $ \frac{c \cdot \pi(y) - z^*}{c \cdot e - z^*} \leq \epsilon  $  within $ \widehat{J}_i $ inner iterations during outer iteration $ i $. 
\vspace{1mm}

To establish (\ref{eqn.cn}), assume $ \widehat{J}_i <  J_i $ and yet the inequality on the right of (\ref{eqn.cn})  does not hold.  (We obtain a contradiction.) For every $ j \leq \widehat{J}_i $, Proposition~\ref{prop.be}  then implies 
\begin{align*}
 \lambda_i^* -  \lambda_{i,j} & > \epsilon \, \,  \frac{c \cdot e - z_i}{c \cdot e - z^*}  \\
                 & =  (1 - \lambda_i^*) \, \epsilon  \quad \textrm{(by (\ref{eqn.ci}))} \; ,
\end{align*}                  
 and hence, using (\ref{eqn.cl}),
 \[  \mathrm{dist}_{z_i}(y_{i,j+1})^2 <  \mathrm{dist}_{z_i}(y_{i,j})^2   - ( \smfrac{3}{4}  - \lambda_i^*) \,  \left( \epsilon / M \right)^2 \; , \]
from which inductively follows
\begin{align*}
   \mathrm{dist}( y_{i,\widehat{J}+1})^2 & <  \mathrm{dist}( y_{i,0})^2 - (\widehat{J}_i +1) \, (\smfrac{3}{4} - \lambda_i^* ) \, ( \epsilon/M)^2 \\
& \leq \mathrm{Dist}^2  - (\widehat{J}_i +1) \, (\smfrac{3}{4} - \lambda_i^* ) \, ( \epsilon/M)^2 \\
& <   0 \; , 
\end{align*}
a contradiction. The claim is established. \vspace{1mm}

Assume $ \epsilon < 1/2 $, the case remaining to be considered. \vspace{1mm}

As each outer iteration $ i $ satisfying $ \lambda_i^* \geq 1/2 $ has only finitely many inner iterations, there must be at least one outer iteration $ i $ satisfying $ \lambda_i^* < 1/2 $.  
Let $ i $ be the first outer iteration for which $ \lambda_i^* < 1/2 $.  From (\ref{eqn.ch})  and (\ref{eqn.cm}), the total number of inner iterations made before reaching outer iteration $ i $ is at most
\begin{equation}  \label{eqn.co} 
\frac{8 (M \, \mathrm{Dist})^2}{ \epsilon } \, \,  \log_{4/3} \left( \frac{1}{ 1 - \frac{z_1 - z^*}{c \cdot e - z^*}} \right)   \; .    
  \end{equation} 

According to (\ref{eqn.cn}), during outer iteration $ i $, the algorithm either achieves its goal within $ \widehat{J}_i $ inner iterations, or the algorithm makes no more than $ \widehat{J}_i $ inner iterations before starting a new outer iteration.  Assume the latter case. Then, for outer iteration $ i + 1 $, the algorithm either achieves its goal within $ \widehat{J}_{i+1} $ inner iterations, or the algorithm makes no more than $ \widehat{J}_{i+1} $ inner iterations before starting a new outer iteration. Assume the latter case. In iteration $ i + 2 $, the algorithm definitely achieves its goal within $ \widehat{J}_{i+2} $ inner iterations, because there cannot be a subsequent outer iteration due, by (\ref{eqn.cj}), to
\[   \smfrac{4}{3} (1 - \lambda_{i+2}^*) \geq \left( \smfrac{4}{3} \right)^3 (1 - \lambda_i^*) >  \left( \smfrac{4}{3} \right)^3 \smfrac{1}{2} > 1 \; . \]

The total number of inner iterations made before the algorithm achieves its goal is thus bounded by the sum of the quantity (\ref{eqn.co})   and 
\begin{align*}
          &    \widehat{J}_i + \widehat{J}_{i+1} + \widehat{J}_{i+2} \\
&  <  \left(  \frac{1}{(1-\frac{1}{2}) - \frac{1}{4}  } + \frac{1}{\frac{4}{3} \, (1 - \smfrac{1}{2}) - \smfrac{1}{4}  } + \frac{1}{\frac{4}{3} \,  \frac{4}{3} \, (1 - \smfrac{1}{2}) - \smfrac{1}{4}}  \right) \, \, \left( \frac{M \, \mathrm{Dist}}{ \epsilon } \right)^2 \\
& \qquad   \textrm{(using $ \smfrac{3}{4} - \lambda_j^* = (1 - \lambda_j^*) - \smfrac{1}{4} $)} \\ 
& < 8 \,  \left( \frac{M \, \mathrm{Dist}}{ \epsilon } \right)^2 \; , 
\end{align*}
completing the proof of the theorem. \hfill $ \Box $

\section{{\bf  General Convex Optimization}} \label{sect.d}

We close with a general example meant to illustrate the flexibility of the preceding development, an example which also serves to highlight a few key differences between our approach and much of the literature on subgradient methods.  \vspace{1mm} 

Let $ f: {\mathcal E} \rightarrow (-\infty , \infty] $ be an extended-valued and lower-semicontinuous convex function.  Consider an optimization problem
\begin{equation}   \label{eqn.da}
  \begin{array}{rl}  \min & f(x) \\
                           \textrm{s.t.} & x \in \mathrm{Feas} \; , \end{array} 
                           \end{equation}  
where $ \mathrm{Feas} = \{ x \in S: Ax = b \} $ and $ S $ is a closed, convex set with nonempty interior.  Assume $ f^* $ -- the optimal value --  is finite. \vspace{1mm}

We explore the complexity ramifications of solving (\ref{eqn.da})  by converting it to an equivalent conic optimization problem and then applying Algorithm 2 of the preceding section (consideration of Algorithm 1 is exactly similar).  \vspace{1mm}

Assume $ \tilde{x} $ is a known point lying in the interiors of both $ S $ and $ \textrm{eff\_dom}(f) $, the effective domain of $ f $ (i.e., where $ f $ is finite).  For convenience of exposition, assume $ \tilde{x} $ is not optimal for (\ref{eqn.da}).  \vspace{1mm}

Assume $ {\mathcal E} $ is endowed with an inner product $ \langle \; , \; \rangle $, and assume $ \| \, \, \| $, the associated norm, satisfies 
\begin{equation}  \label{eqn.db}
      \{ x : \| x - \tilde{x} \| \leq 1 \textrm{ and }    Ax = b \} \subseteq \mathrm{Feas}  \cap \mathrm{eff\_dom (f)} \; .
      \end{equation} 
Let $ \hat{f}    $ be a known scalar satisfying $ \hat{f}   \geq f(x) $ for all $ x $ in the set on the left of (\ref{eqn.db}). (The value $ \hat{f} $ is required to be known because it together with $ \tilde{x} $  will determine the distinguished direction $ e $.)\vspace{1mm}

Let $ D $ denote the diameter of the sublevel set $ \{ x \in \mathrm{Feas}: f(x) \leq f( \tilde{x}) \} $.  Assume $ D $ is finite, implying the optimal value $ f^* $ of (\ref{eqn.da}) to be attained by some feasible point.  \vspace{1mm}

 For later reference, observe that assumption (\ref{eqn.db})  and the convexity of $ f $ imply
\[     f( \tilde{x})  \leq \smfrac{1}{D +1} f^* + \smfrac{D }{D +1} \hat{f}   \; , \]
which in turn implies
\begin{equation}  \label{eqn.dc} 
\frac{1}{1 - \frac{f( \tilde{x}) - f^*}{ \hat{f} - f^*}} \, \leq \, D+1 \; . 
\end{equation}
 
As $ S $ is closed and convex, there exists a closed, convex cone $ {\mathcal K}_1 \subseteq {\mathcal E} \times \mathbb{R} $ for which $ S = \{ x: (x,1) \in {\mathcal K}_1 \} $. \vspace{1mm}  

For an extended-valued function to be lower semicontinuous is equivalent to its epigraph being closed.  Since the epigraph for $ f $ is convex, there thus exists a closed, convex cone $ {\mathcal K}_2 \subseteq {\mathcal E} \times \mathbb{R} \times \mathbb{R} $ for which
\[   \mathrm{epi}(f) := \{ (x,t): f(x) \leq t \} = \{ (x,t): (x,1,t) \in {\mathcal K}_2 \} \; . \]
Note 
\begin{equation}  \label{eqn.dd}
  t > f(\tilde{x}) \quad \Rightarrow \quad  (\tilde{x}, 1,t) \in \int({\mathcal K}_2) \; , 
  \end{equation}
a consequence of the assumption $ \tilde{x} \in \int(\mathrm{eff\_dom (f)}) $. \vspace{1mm}

Let 
\[  {\mathcal K} := \{ (x,s,t): (x,s) \in {\mathcal K}_1 \textrm{ and } (x,s,t) \in {\mathcal K}_2 \} \; . \]
Clearly, the optimization problem (\ref{eqn.da})  is equivalent to
\begin{equation}  \label{eqn.de}
  \begin{array}{rl}
 \min_{x,s,t} & t \\
      \textrm{s.t.} & Ax = b \\
                    & s = 1 \\
   & (x,s,t) \in {\mathcal K} \; , \end{array} \end{equation} 
and has the same optimal value, $ f^* $.  The conic program (\ref{eqn.de})  is of the same form as CP, the focus of preceding sections. \vspace{1mm}
  
Observe for all scalars $ t $,
\begin{equation}  \label{eqn.df}
 \mathrm{Level}_{t} = \{ (x,1, t ): x \in \mathrm{Feas} \textrm{ and }   f(x) \leq t   \} \; . \end{equation}
     
To apply Algorithm 2, a distinguished direction $ e $ is needed, and a computational inner product should be specified. Additionally, an input $ \bar{x} $ is required.    \vspace{1mm}

Let $ e = ( \tilde{x}, 1, \hat{f} \, ) $, which lies in the interior of $ {\mathcal K} $, due to (\ref{eqn.dd})  and $ \tilde{x} \in \int(S) $.  This distinguished direction, along with the cone $ {\mathcal K} $, determine the map $ (x,s,t) \mapsto \lambda_{\min}(x,s,t) $ on $ {\mathcal E} \times \mathbb{R} \times \mathbb{R} $. \vspace{1mm}          

Choose the computational inner product on $ {\mathcal E} \times \mathbb{R} \times \mathbb{R} $ to be any inner product that assigns to pairs $ (x_1,0,0), (x_2,0,0) $ the value $ \lin x_1, x_2 \rin $ (the original inner product on $ {\mathcal E} $). By (\ref{eqn.db})  and (\ref{eqn.df}), 
 the level set containing $ e $ then satisfies 
\[  \mathrm{Level}_{\hat{f}} \cap B(e,1) \subseteq {\mathcal K} \; , \]
and hence, by Proposition~\ref{prop.ca}, the Lipschitz constant is at most 1 for the map $ (x,s,t) \mapsto \lambda_{\min}(x,s,t) $ restricted to $ \mathrm{Affine}_t $, for every $ t $.  \vspace{1mm}

Choose the input $ \bar{x}  $ to Algorithm 2 as $ \bar{x}  = ( \tilde{x}, 1, f( \tilde{x})) $, which clearly is feasible for the conic program (\ref{eqn.de}).  Note that (\ref{eqn.df})  then implies the horizontal diameter of the relevant sublevel set for the conic program satisfies 
\[  \mathrm{Diam}_{f(\tilde{x})} =  D   \]
(i.e., is equal to the diameter of the sublevel set $ \{ x \in \mathrm{Feas}: f(x) \leq f( \tilde{x}) \} $). \vspace{1mm}

Applying Algorithm 2 results in a sequence of iterates $ (x_k, 1, t_k) $ for which the projections $ (x_k',1,t_k') := \pi(x_k,1,t_k) $ satisfy $ x_k' \in \mathrm{Feas} $ and $ f(x_k') \leq t_k' $ (simply because $ (x_k',1,t_k') $ is feasible for the conic program (\ref{eqn.de})). \vspace{1mm}

   Since $ (x_0,1,t_0) = (\tilde{x},1,f(\tilde{x})) $, we have $ (x_0',1,t_0') = (x_0,1,t_0) $ -- in particular, we have $ t_0' = f( \tilde{x}) $.  Consequently, the sequence of points $ x_k' $ not only lie in $ \mathrm{Feas} $, but by Theorem~\ref{thm.ch}  satisfy   
\begin{align}  
 & \ell  \geq 8 \, D^2 \, \left( \, \frac{1}{\epsilon^2} \, + \, \frac{1}{\epsilon} \, \log_{4/3} \big( (D+1) \, (1 - \epsilon ) \big) \,  + 1 \, \right)   \label{eqn.dg}   \\
&  \quad                \quad \Rightarrow \quad \min_{k \leq \ell } \,   \frac{f(x_k') - f^*}{\hat{f} - f^*} \, \leq \, \epsilon \; ,   \label{eqn.dh} 
 \end{align}
where for (\ref{eqn.dg})  we have used (\ref{eqn.dc}), and for (\ref{eqn.dh})  have used $ f(x_k') \leq t_k' $. \vspace{1mm}

Deserving of emphasis is that the only projections made are onto the subspace $ 
  \{ (x,s,t): Ax = 0 \textrm{ and } s = 0 \} $ (the same subspace at every iteration, a situation for which preprocessing is effective, especially if the number of equations is relatively small). This differs from much of the subgradient method literature where, for example, commonly required is projection onto $ \mathrm{Feas} $ for each iterate landing outside $ \mathrm{Feas} $. (A projection swamps the cost of an iteration except for especially simple sets $ \mathrm{Feas} $, such as a box, a ball, or an affine space.)    
  \vspace{1mm}

Another difference, deserving perhaps of even more emphasis, is that the bound (\ref{eqn.dg})  is independent of a Lipschitz constant for $ f $. In fact, no Lipschitz constant is implied by the assumptions, as is seen  by considering the family of univariate cases in which $ \mathrm{Feas} = \mathbb{R} $, $ A = b = 0 $, and $ f $ is allowed to be {\em any}  lower-semicontinuous convex function with 
\[  [0,2] \subseteq \mathrm{eff\_dom}(f) \subseteq [0, \infty) \; ,  \]
and which is strictly increasing at $ 0 $. The assumptions are fulfilled by choosing $ \tilde{x} = 1 $ and $ \hat{f} = f(2) $, and by using the standard inner product (i.e., multiplication), in which case $ D = 1 $.    The bound (\ref{eqn.dg})  is then $ \ell  \geq 8 ( \smfrac{1}{\epsilon^2} + \smfrac{1}{\epsilon }  \log_{3/2}(2(1 - \epsilon ) ) + 1 ) $, whereas the error bound (\ref{eqn.dh})  is $ \frac{f(x_k') - f(0)}{f(2) - f(0)} \leq \epsilon $.  Clearly, even for the restriction of $ f $ to $ [0,2] $, nothing is implied about the Lipschitz constant other than trivial {\em lower} bounds such as $ L \geq \smfrac{1}{2} (f(2) - f(0))\; . $   \vspace{1mm}

For the approach presented herein, Lipschitz continuity matters only with regards to the function $ (x,s,t) \mapsto \lambda_{\min}(x,s,t) $ restricted to $ \mathrm{Affine}_{t} $, which is guaranteed to have Lipschitz constant at most 1 (due to assumption (\ref{eqn.db})).   On the other hand, the error bound (\ref{eqn.dh})  is measured relatively, whereas in traditional subgradient-method literature relying on a Lipschitz constant for the objective function $ f $, error is specified absolutely.

\bibliographystyle{plain}
\bibliography{subgrad_framework}

  \end{document}